\documentclass[11pt,twoside]{amsart}
\usepackage{pgfplots}
\pgfplotsset{compat=1.15}
\usepackage{mathrsfs}
\usetikzlibrary{arrows}
\usepackage{amsmath, amssymb, amsfonts, amsthm,latexsym,,mathtools,,hyperref,fancyhdr,comment, tabularht, tabularx, longtable, mathrsfs, cleveref,booktabs, color, geometry,bm}
\usepackage{enumitem}
\usepackage{multicol}
\usepackage{hyperref,url}
\usepackage{tikz}
\tikzstyle{vertex}=[circle, draw, inner sep=0pt, minimum size=6pt,fill=black]
\newcommand{\vertex}{\node[vertex]}
\long\def\symbolfootnote[#1]#2{\begingroup%
	\def\thefootnote{\fnsymbol{footnote}}\footnote[#1]{#2}\endgroup}

\allowdisplaybreaks
\def \N {{\mathbb{N}}}
\def \Z {{\mathbb{Z}}}

\def \F {{\mathbb{F}}}

\def \Ga {{\Gamma }}

\newtheorem*{theorem*}{Theorem}
\newtheorem{theorem}{Theorem}[section]
\newtheorem{cor}[theorem]{Corollary}

\newtheorem{lemma}[theorem]{Lemma}
\newtheorem{ex}[theorem]{Example}

\newtheorem*{ex*}{Example}

\newtheorem{pro}[theorem]{Proposition}
\newtheorem{problem}[theorem]{Problem}
\date{}

\begin{document}
\title{Cyclic Subgroup Graph of a Group}
\author{Peter J. Cameron}
	\address{School of Mathematics and Statistics, University of St Andrews, Fife KY16 9SS, U.K.}
	\email{pjc20@st-andrews.ac.uk}
	\author{A. satyanarayana Reddy}
	\address{Shiv Nadar Institution of Eminence, NH-91, Dadri, Gautam Buddha Nagar}
	\email{satya.a@snu.edu.in}
    
	\author{Khyati Sharma}
	\address{Shiv Nadar Institution of Eminence, NH-91, Dadri, Gautam Buddha Nagar}
	\email{khyatisharma0907@gmail.com}

	\today 

\begin{abstract}
   A {\em cyclic subgroup} graph of a group is a graph whose vertices are cyclic subgroups of the group, and two distinct vertices are adjacent if one is a maximal subgroup of the other. In this paper, we examine several graph-theoretic properties of cyclic subgroup graphs, including bipartiteness, regularity, girth, diameter, and vertex degrees. We also classify all groups whose cyclic subgroup graphs are connected, complete, star graphs, trees, or paths. Further, we examine the cyclic subgroup graph for specific families of groups such as cyclic, dihedral, dicyclic, generalized quaternion, minimal non-cyclic and nilpotent groups. We also discuss the relationship of cyclic subgroup graph with the power graph and the enhanced power graph. Additionally, we give a relationship between a cyclic subgroup graph and a poset.
\end{abstract}
\maketitle
\section{Introduction}\label{sec1}
Graphs associated with algebraic structures are
an effective tool for relating the algebraic structure to its combinatorial and geometric properties. These properties have valuable applications in the applied sciences. Additionally, studying these graphs helps us to understand the properties of the associated algebraic structures. There are various graphs arising from group-theoretic properties that have been studied in the literature, such as commuting graphs~\cite{akbari2004commuting}, power graphs~\cite{Chakrabarty2009,kumar2021recent}, enhanced power graphs~\cite{AalipourEtAl2017,bera2022proper,bera2021connectivity}, and comaximal subgroup graphs~\cite{Akbari_Miraftab_Nikandish_2017,das2024co,das2025co}. An extensive survey on different graphs defined on groups has been done by P. J. Cameron~\cite{PaterCameron}.

One such graph is known as the {\em cyclic subgroup graph} $\Gamma (G)$ of a group $G$. Its vertex set is $C(G)$, the collection of all cyclic subgroups of $G$, and two vertices $H_1$ and $H_2$ are adjacent if $H_1\leq H_2$, and there is no subgroup $K$ such that $H_1< K< H_2$.

The study of cyclic subgroup graphs was initiated by M. Tărnăuceanu in~\cite{tuarnuauceanu2023number}, where an explicit formula for determining the number of edges was introduced. By definition, it is easy to see that these graphs are induced subgraphs of the {\em subgroup graph} $L(G)$ of $G$, that is, the graph whose vertices are the subgroups of $G$, two vertices $H_1$ and $H_2$ adjacent if $H_1\leq H_2$ and there is no subgroup $K$ such that $H_1<K<H_2$. For more details about the subgroup graph, we refer~\cite{bohanon2006finite, lucchini2021genus, schmidt2006planar}.

Cyclic subgroup graphs can be defined for infinite groups as well, but in this work, we restrict ourselves to finite groups only. Here, all the graphs are undirected without multiple edges and loops. The organization of the article is as follows. In the rest of this section, we provide a few notations, preliminary definitions, and results that we use throughout this work. In Section~\ref{sec:2}, we first characterize when the cyclic subgroup graph $\Gamma (G)$ is bipartite, connected, complete, a tree, regular, and a cograph, etc. We also study various graph-theoretic parameters like diameter, degree, and girth for cyclic, dihedral, dicyclic, and generalized quaternion groups. Section~\ref{sec:nil} studies properties of the cyclic subgroup graph for nilpotent groups. A non-cyclic group is {\em minimal non-cyclic} if all its proper subgroups are cyclic. Section~\ref{sec:3} explores properties of the cyclic subgroup graph associated with minimal non-cyclic groups. In Section~\ref{sec:relationship}, the relationship of the cyclic subgroup graph with the power graph and the enhanced power graph is given. At last, in Section~\ref{sec:poset} we discuss the relationship of the cyclic subgroup graph with the poset of cyclic subgroups ordered by inclusion.
\subsection{Notations and Preliminaries}
The notations $\Z_n$, $D_{2n}$, $Q_{2^n}$, and $Dic_{n}$ denote the cyclic group of order $n$, dihedral group of order $2n$, generalized quaternion group of order $2^n$, and dicyclic group of order $4n$, respectively. A group $G$ is {\em CLT} if it has a subgroup corresponding to every divisor of the order of $G$. The {\em diameter} of a connected graph $\Gamma$ is the maximum distance between any two vertices in $\Gamma$, we denote it by $diam (\Gamma)$. The {\em degree} of a vertex $v$ in a graph $\Gamma$ is the number of edges incident to the vertex $v$, and it is denoted by $deg_{\Gamma}(v)$. If the graph $\Gamma$ is known, then we can simply write it as $deg (v)$. The {\em maximum} and {\em minimum degree} in the graph $\Gamma$ is denoted by $\Delta (\Gamma)$ and $\delta (\Gamma)$, respectively. A vertex of a graph is {\em pendant} if its degree is one. The {\em girth} $(\Gamma)$ denotes the length of the smallest cycle in $\Gamma$, if a graph does not contain a cycle, then the girth is infinite. A connected graph $\Gamma$ is {\em Eulerian} if there exists a closed trail containing every edge of $\Gamma$. A graph $\Gamma$ is {\em regular} if the degree of each vertex is the same. A {\em cograph} is a type of graph that contains no induced path on four vertices. A subset of vertices in a graph where no two vertices are adjacent is known as an {\em independent set}, and the maximum cardinality of an independent set is the {\em independence number}. A subset $S$ of vertices of a graph $\Gamma$ is called a {\em dominating set} if any vertex of $\Gamma$ is either in $S$ or adjacent to some vertex in $S$. The minimum cardinality of a dominating set is called {\em dominating number} $\gamma$ of $\Gamma$. For other basic definitions used in this article, we refer to the book~\cite{west2001introduction}. The numbers $p_1,p_2,\ldots$ denote distinct prime numbers. In this paper, the exponent of a group $G$ is defined as the smallest positive integer $n$ such that $g^n=e$ for every $g\in G$. On the other hand, the exponent of a prime $p$ in a natural number $m$ refers to the multiplicity of $p$ in the prime factorization of $m$. The following results are used throughout this paper.
\begin{theorem}\label{maximal}[{M. Hall, \cite{hall2018theory}}]
 Every maximal subgroup of a nilpotent group is normal with a prime index.   
\end{theorem}
\begin{theorem}\label{sylowspecial}[{Frobenius, \cite{brown1988generalization}}] If $p$ is a prime number and $p$ divides the order of a group $G$, and $n_p$ denotes the number of subgroups of order $p$ in $G$, then $n_p\equiv 1 \mod p$.
\end{theorem}
\section{Properties of $\Ga(G)$}\label{sec:2}
In this section, we discuss various properties of $\Ga (G)$. We start with an example of $\Ga (G)$.
\begin{ex}\label{example1}
Consider the group $G\cong D_{12}=\langle r, s:r^6=s^2=e, srs=r^{-1}\rangle.$ The set of all cyclic subgroups of $G$ is $C(G)=\{\{e\}, \Z_2=\langle s\rangle, \Z_2=\langle sr\rangle, \Z_2=\langle sr^2\rangle, \Z_2=\langle sr^3\rangle, \Z_2=\langle r^3\rangle, \Z_2=\langle sr^4\rangle, \Z_2=\langle sr^5\rangle, \Z_3=\langle r^2\rangle, \Z_6=\langle r\rangle\}$, and the corresponding cyclic subgroup graph $\Gamma (G)$ is given in Figure~\ref{fig45}.
\begin{figure}[!htbt]
  \centering
  \begin{tikzpicture}
  [scale=0.80]
	\vertex (a) at (0,0)[label=below:$\{e\}$]{};
	\vertex (b) at (-2,2)[label=above:$\Z_2$] {};
	\vertex (c) at (2,2) [label=above:$\Z_{3}$] {};
	\vertex (d) at (0,4)[label=above:$\Z_{6}$]  {};
	\vertex (e) at (-4,1)[label=above:$\Z_{2}$]  {};
        \vertex (f) at (-3,1)[label=above:$\Z_{2}$]  {};
	\vertex (g) at (-2,1)[label=above:$\Z_{2}$]  {};
        \vertex (h) at (2,1)[label=above:$\Z_2$]  {};
        \vertex (i) at (3,1)[label=above:$\Z_{2}$]  {};
        \vertex (j) at (4,1)[label=above:$\Z_{2}$]  {};
	\path
	  	(a) edge(b)
		(a) edge(c)
		(b) edge(d)
            (c) edge(d)
            (a) edge(e)
		(a) edge(f)
            (a) edge(g)
		(a) edge(h)
		(a) edge(i)
		(a) edge(j)

	 ;
\end{tikzpicture}
\caption{Cyclic subgroup graph of $D_{12.}$}
\label{fig45}
\end{figure}
\end{ex}
 In  Table~\ref{tab:vertexedges}, we list the number of vertices and edges of cyclic subgroup graphs for some particular families of groups by using~\cite{tuarnuauceanu2023number} for a better understanding of the structure of these graphs.
\begin{table}[ht!]
\centering
\begin{tabular}{|c|c|c|}
 \hline
 $\mathbf{G}$ & \textbf{Number of vertices} & \textbf{Number of edges}\\
 \hline
  & &\\
  $\Z_n$, $n=p_1^{a_1}p_2^{a_2}\cdots p_k^{a_k}$ &
   $\prod_{i=1}^{k}(a_i+1)$ & $\bigg(\sum_{i=1}^{k}\frac{a_i}{a_i+1}\bigg)\prod_{i=1}^{k}(a_i+1)$\\
   \hline
    & &\\
  $D_{2n}$, $n=p_1^{a_1}p_2^{a_2}\cdots p_k^{a_k}$ & $\prod_{i=1}^{k}(a_i+1)+n$ &
    $\bigg(\sum_{i=1}^{k}\frac{a_i}{a_i+1}\bigg)\prod_{i=1}^{k}(a_i+1)+n$\\
    \hline
     & &\\
  $Q_{2^n}$, $n\geq 3$ & $2^{n-2}+n$ & $2^{n-2}+n-1$\\
  \hline
   & &\\
  $Dic_{n}$, $2n=2^{a_1}p_2^{a_2}\cdots p_k^{a_k}$ & $\prod_{i=1}^{k}(a_i+1)+n$ & $\bigg(\sum_{i=1}^{k}\frac{a_i}{a_i+1}\bigg)\prod_{i=1}^{k}(a_i+1)+n$\\
  \hline
   & &\\
  $\Z_{p^a}\times \Z_q\times \Z_q$ & $(a+1)(q+2)$ & $(a-1)(2q+3)+(3q+4)$\\

  \hline
 \end{tabular}\\
 \caption{}
           \label{tab:vertexedges}
       \end{table}
\begin{pro}\label{bipartite}
           Cyclic subgroup graphs are always bipartite.
       \end{pro}

       \begin{proof}
This proof introduces a concept we will use again later. A \emph{hommorphism}
from a graph $\Gamma_1$ to a graph $\Gamma_2$ is a map from the vertex set of
$\Gamma_1$ to that of $\Gamma_2$ which carries edges to edges. Now a graph is
bipartite if and only if it has a hommomorphism onto $K_2$, the graph with a
single edge: the bipartition consists of the inverse images of the two vertices
of the edge.

Let $G$ be a finite group of order $n=p_1^{a_1}p_2^{a_2}\dots p_k^{a_k}$, where $a_i\geq 1$, and let $\Ga(G)= (V, E)$ be the corresponding cyclic subgroup graph of $G$. If $H$ is a cyclic subgroup of $G$ of order $p_1^{b_1}p_2^{b_2}\dots p_k^{b_k}$, then we define $e(H)=\sum_{i=1}^{k}b_i$. 

Now let $e$ and $o$ label the vertices of a complete graph $K_2$, and map the
subgroup $H$ to $e$ if $e(H)$ is even, and to $o$ if it is odd. It is easy to
see that this is a graph homomorphism. Hence $\Ga(G)$ is bipartite.
\end{proof}

In fact we can always find a homomorphism to a larger graph. Recall that the
\emph{Cartesian product} of graphs $\Gamma_1,\ldots,\Gamma_k$ is the graph
whose vertex set is the Cartesian product of the vertex sets of these graphs,
with vertices $(v_1,\ldots,v_k)$ and $(w_1,\ldots,w_k)$ adjacent if $v_i=w_i$
for all but exactly one value of $i\in\{1,\ldots,k\}$.

\begin{pro}
\begin{enumerate}
\item Let $n=p_1^{a_1}\cdots p_k^{a_k}$, where $p_1,\ldots,p_k$ are distinct
primes and the exponents $a_1,\ldots,a_k$ are positive.
Then the cyclic subgroup graph of $\Z_n$ is
isomorphic to the Cartesian product of paths of lengths $a_1,\ldots,a_k$.
\item Let $G$ be a group with exponent $n$ as in the preceding part. Then
the cyclic subgroup graph of $G$ has a homomorphism into the cyclic subgroup
graph of $\Z_n$.
\end{enumerate}
\end{pro}

\begin{proof}
\textbf{(1)} Label the vertices of the $i$th path $0,1,\ldots,a_i$. Then map the
unique subgroup of order $p_1^{b_1}\cdots p_k^{b_k}$
(where $0\le b_i\le a_i$ for all $i$) to the vertex $(b_1,\ldots,b_k)$ of
the Cartesian product of the paths. It is readily checked that this map
is a graph isomorphism.

\textbf{(2)} The homomorphism can be described briefly by saying it maps a cyclic
subgroup of order $m$ in $G$ to the unique cyclic subgroup of order $m$
in $\Z_n$.
\end{proof}


\begin{theorem}\label{path}
    Let $G$ be a finite group. Then
    \begin{enumerate}[label={(\arabic*)}]
    \item\label{Thm:part:1} $\Ga(G)$ is connected and it is a complete graph if and only if $G$ is either trivial or $G\cong \Z_p.$
        \item\label{Thm:part:2} $\Ga (G)$ is a path graph if and only if $G\cong \Z_{p^k}$, where $k\in \N.$
        \item\label{Thm:part:3} $\Gamma (G)$ is a cycle graph if and only if $G\cong \Z_{pq}$.
        \item\label{Thm:part:4} $\Gamma (G)$ is a star graph if and only if 
$G$ is isomorphic to $Z_{p^2}$, $Q_8$, or a group in which every element has prime order.
        \end{enumerate}
        \end{theorem}
        \begin{proof}
     \textbf{\ref{Thm:part:1}}   
     It is easy to see that there is a path from $\{e\}$ to $H$ as $\{e\}\sim H_1\sim H_2 \sim \dots \sim H_{m-1}\sim H_m=H,$ where $|H_{i+1}/H_i|$ is a prime divisor of $|H|.$ Hence $\Gamma(G)$ is connected. Since $\Ga(G)$ is bipartite and connected, $\Ga(G)$ is a complete graph if and only if $|C(G)|\le 2$. Hence $G$ is either trivial or $\Z_p$.

     \noindent
\textbf{\ref{Thm:part:2}} If $G\cong \Z_{p^k}$, where $k\in \N,$ then it is easy to see that $\Ga(G)$ is a path graph. Conversely, suppose that $\Ga(G)$ is a path graph. First, we claim that the degree of vertex $\{e\}$ is $1.$ If the  degree of vertex $\{e\}$ is not $1,$ then it is $2$ as  $\Ga(G)$ is a path graph. Consequently $G$ has  exactly two  subgroups $H$ and $K$ of order $p$ and $q$ respectively, where $p$ and $q$ distinct primes. Hence the vertices $\{\{e\},H,K,HK\}$ form a cycle graph $C_4.$ This implies that the degree of $\{e\}$ is $1$ and $|G|=p^n$ for some prime $p$. Now, our goal is to show that $G$ is a cyclic group. Let us assume that $G$ is not cyclic. Since $\Gamma (G)$ is a path, then the length of the path is at most $n$, and the number of cyclic subgroups of $G$ is at most $n+1$. This is a contradiction as by Richard's theorem~\cite{richards1984remark}, the number of cyclic subgroups of $G$ is more than $\tau (n)=n+1$, where $\tau (n)$ denotes the number of divisors of $n$. Hence $G$ is a cyclic $p$-group.

          \noindent \textbf{\ref{Thm:part:3}} If $G=\Z_{pq},$ then $C(G)=\{\{e\},\Z_p,\Z_q, \Z_{pq}\}$ and $\Ga(G)$ is the cycle graph $C_4.$ Conversly if $\Ga(G)$ is a cycle graph, then degree of $\{e\}$ is $2$ and from the proof of part-\ref{Thm:part:1},  $\Ga(G)$ contains cycle $C_4.$ Hence the result follows.

     \noindent \textbf{\ref{Thm:part:4}} In a star graph, all the vertices are of degree one except the central vertex. Now, there are two cases.
     
    \textbf{Case \;1:} If $\{e\}$ is the central vertex, then $\Gamma (G)$ is the star graph if all the non-trivial cyclic subgroups of $G$ are of prime order, otherwise there exists at least one more vertex of degree more than $1$. A complete classification of these groups is given in \cite{Cheng1993}.
    
\textbf{Case \;2:} If $\{e\}$ is not the central vertex, then the degree of $\{e\}$ is $1$, which implies that $G$ is a $p$-group with unique prime ordered subgroup say $\Z_p$. In this case, the vertex $\Z_p$ is the central vertex as either $G\cong \Z_p$ or degree of $\Z_p$ will be at least $2$. Since $G$ has a unique prime ordered subgroup, then the Sylow $p$-subgroup of $G$ is either cyclic or a generalized quaternion \cite[Theorems 8.5,8.6]{burnside_2012}. If Sylow $p$-subgroup of $G$ is cyclic, then it is easy to see that $G\cong \Z_{p^2}$. If Sylow $p$-subgroup is generalized quaternion, then $G\cong Q_8$ otherwise, $\Gamma (G)$ has at least two vertices of degree more than two. 
\end{proof}

Now we show the existence of induced paths and cycles in $\Ga(G)$.

\begin{pro}
Let $G$ be a finite group containing an element of order
$n=p_1^{a_1}\dots p_k^{a_k}$, where $p_1,\ldots,p_k$ are distinct primes
and $a_1,\ldots,a_k>0$. Suppose that either $k\ge2$, or $a_i>1$ for some $i$.
Then
\begin{enumerate}
\item $\Ga(G)$ contains an induced cycle of length
$2r=2(a_1+\cdots+a_k)$;
\item If $G$ contains two distinct subgroups of order $n$, then $\Ga(G)$
contains an induced path of length $2r$.
\end{enumerate}
\end{pro}

\begin{proof}
The sequence of subgroups of $\langle g\rangle$ of orders
$1,p_1,\ldots,p_1^{a_1},p_1^{a_1}p_2,\ldots,p_1^{a_1}p_2^{a_2},\cdots,n$
is an induced path of length $a_1+\cdots+a_k$ from the identity subgroup
to $G$.

If $k>2$, then another induced path is obtained by reversing the order of
the primes, that is, taking subgroups of orders
$1,p_k,\ldots,p_k^{a_k},p_k^{a_k}p_{k-1},\ldots,n$. It is readily checked
that no edges join vertices in one path to those of the other except for
the end points; so we have an induced cycle of length $r$.

If $k=2$ and (without loss) $a_1>1$, the same construction works.

Now suppose that there is another cyclic subgroup $\langle h\rangle$ of
order~$n$. Then take the first path from the trivial group to
$\langle g\rangle$, and the second path from trivial group to 
$\langle h\rangle$; the induced subgraph is a path of length $r$.
\end{proof}

An {\em EPPO group} is a group all of whose elements have prime power order. These groups are very well studied in the literature. G. Higman studied solvable EPPO groups~\cite{Higman1957EPPO} in 1957. Then, in 1961, M. Suzuki~\cite{Suzuki1961NilpotentCentralizers}
classified the simple EPPO groups. Later, in 1981, R. Brandl classified all EPPO groups in~\cite{Brandl1981PrimePowerOrder}. We note that the list includes a few finite simple groups: $A_5$, $A_6$, $\mathrm{PSL}(2,7)$, $\mathrm{PSL}(2,8)$,
$\mathrm{PSL}(2,17)$, $\mathrm{Sz}(8)$, $\mathrm{Sz}(32)$, and 
$\mathrm{PSL}(3,4)$.

\begin{theorem}
    The cyclic subgroup graph $\Gamma (G)$ of a group $G$ is a tree if and only if $G$ is an EPPO group.
\end{theorem}

\begin{proof}
    If $G$ is not an EPPO group, then it contains an element with order divisible by at least two different primes $p$ and $q$, hence $G$ has an element of order $pq$, which gives a cycle in the cyclic subgroup graph $\Gamma (G)$. Conversely, if the graph contains a cycle, choose a subgroup of maximum order from the vertices of the cycle. This subgroup must cover at least two subgroups, because of the cycle and by the definition of $\Gamma (G)$  and so its order cannot be a prime power. Hence, the result holds.

\end{proof}
\noindent The following results computes degree of each vertex of $\Ga(G)$, where $G\in \{\Z_n, D_{2n}, Q_{2^n}, Dic_n\}.$

\begin{theorem}\label{degree}
        Let $G\in \{\Z_n, D_{2n}, Q_{2^n}, Dic_n\}$, where $n=p_1^{a_1}p_2^{a_2}\dots p_k^{a_k}$, and $a_i\geq 1$. Also let $H$ be a cyclic subgroup of $G$ such that
       $|H|=p_1^{b_1}p_2^{b_2}\dots p_k^{b_k}$, and  $m(H)=|\{i|b_i=0\}|$ and $r(H)=|\{i|b_i=a_i\}|$ for all $1\leq i\leq k$. Then, the degree of $H$ is listed in Table~\ref{tab:degree}.
        \begin{table}[ht!]
        \centering
\begin{tabular}{|c|c|}
 \hline
 $\mathbf{G}$ & $\mathbf{deg (H)}$\\
 \hline
  $\Z_n$ & $2k-m(H)-r(H).$\\
 \hline
  $D_{2n}$ &
    $\begin{cases}
            1 & \mbox{if $H$ is generated by a reflection},\\
             k+n & \mbox{if $H=\{e\}$},\\
              2k-m(H)-r(H) & \mbox{if $H\ne \{e\}$ and is generated by a rotation.}\\
              \end{cases}$ \\
              \hline
  $Q_{2^n}$ & $\begin{cases}
            2^{n-2}+1 & \mbox{if $H\cong \Z_2$},\\
             2 & \mbox{if $\{e\}\lneq H\lneq \Z_{2^{n-1}}$ and $H\not \cong \Z_2$},\\
              1 & \mbox{otherwise}.
              \end{cases}$\\
              \hline
  $Dic_{n}$ & $\begin{cases}
            k+n+1 & \mbox{if $H\cong \Z_2$},\\
             2k-m(H)-r(H) & \mbox{if $H\leq \Z_{2n}$, where $n$ is even and  $|H|\ne 2$},\\
               2(k+1)-m(H)-r(H) & \mbox {if $H\leq \Z_{2n}$, where $n$ is odd,}\\
               & |H|=2^{b_{k+1}}p_1^{b_1}p_2^{b_2}\cdots  p_k^{b_k}\ne 2,\\
              1 & \mbox{otherwise}.\\
              \end{cases}$\\

  \hline
 \end{tabular}\\
 \caption{Degree of $\Gamma (G)$}
           \label{tab:degree}
       \end{table}
   \end{theorem}
   \begin{proof}
   We know that in a cyclic group, every maximal subgroup is of prime index. We assume $m=m(H)$ and $r=r(H)$ whenever $H$ is clear from the context.  Now, we have the following cases:\\
       $\mathbf{G\cong \Z_n}$, where $n=p_1^{a_1}p_2^{a_2}\dots p_k^{a_k}$, and $a_i\geq 1$. Let $H$ be a subgroup of $G$ of order $p_1^{b_1}p_2^{b_2}\dots p_k^{b_k}.$ Then, we have the following subcases:
        \begin{enumerate}[label={(\arabic*)}]
            \item If $0<b_i<a_i, \forall 1\leq i \leq k,$ then $H$ is adjacent to the subgroups of order $p_1^{b_1}p_2^{b_2}\dots p_i^{b_i-1}\dots p_k^{b_k}$ and $p_1^{b_1}p_2^{b_2}\dots p_i^{b_i+1}\dots p_k^{b_k}$ for all $1\leq i \leq k.$ That is $H$ adjacent to $2k$ subgroups in $\Gamma (G).$ Therefore the degree of $H$ in $\Gamma (G)$ is $2k.$

            \item If $b_i<a_i, \forall \;1\leq i \leq k$   and assume that $b_1, b_2,\ldots b_m$ are zero, then $H$ is adjacent to  the subgroups with order $p_t{p_{m+1}}^{b_{m+1}}\dots p_k^{b_k}$  
  for all $1\leq t\leq m$ and ${p_{m+1}}^{b_{m+1}}\dots p_i^{b_{i+1}}\dots p_k^{b_k}$ for all ${m+1}\leq i \leq k$, also $H$ is adjacent  to subgroup of order ${p_{m+1}}^{b_{m+1}}\dots p_i^{a_{i-1}}\dots p_k^{b_k}$ for all $m+1\leq i\leq k$. Therefore, the degree of $H$ in $\Gamma (G)$ is $2k-m.$

            \item If $b_i > 0, \forall 1\leq i\leq k$, and assume that $b_i=a_i,\forall \; 1\leq i \leq r$, then $H$ is adjacent to the subgroups of order  $p_1^{b_1}p_2^{b_2}\dots p_i^{b_i-1}\dots p_k^{b_k}$ for all $1\leq i\leq k$, and $p_1^{b_1}\dots p_r^{b_r}{p_{r+1}}^{b_{r+1}}\dots p_i^{b_{i+1}}\dots p_k^{b_k}$ for all $r+1\leq i\leq k$. Therefore, the degree of $H$ in $\Gamma (G)$ is $2k-r$.

            \item If $m$ number of $b_i$'s are zero and $r$ number of $b_i=a_i$, then by using the previous cases the degree of $H$ in $\Ga (G)$ is $2k-m-r$.
        \end{enumerate}
          $\mathbf{G\cong D_{2n}}$, where $n=p_1^{a_1}p_2^{a_2}\dots p_k^{a_k}$, and $a_i\geq 1$ and let $H$ be a non-trivial subgroup of $D_{2n}$ of order $p_1^{b_1}p_2^{b_2}\dots p_k^{b_k}$ generated by a rotation. Then, by using the first part, we see that the degree of $H$ in $\Gamma (D_{2n})$ is $2k-m-r.$  If $H$ is generated by a reflection, then $|H|=2$, and $H$ is not contained in any other cyclic subgroup. Therefore, these vertices are pendant in $\Ga(D_{2n}).$ Also, $D_{2n}$ has exactly $k+n$ subgroups of prime order, and $\{e\}$ is a maximal subgroup of all these subgroups, thus, the degree of $\{e\}$ in $\Gamma (G)$ is $k+n$.\\
          
          $\mathbf{G\cong Q_{2^n}}$, where $n\geq 4$. By Theorem~$4.2$ of \cite{conrad2014generalized}, $Q_{2^n}$ has a unique cyclic subgroup of order $2^{n-1}$ and every element outside that has order $4$. Then, it is easy to see that $Q_{2^n}$ has $2^{n-2}+1$ cyclic subgroups of order $4$. Therefore degree of $\Z_2$ in $\Gamma (G)$ is $2^{n-2}+2$, and degree of $H$ is $2$ if and only if $\{e\}\lneq H\lneq \Z_{2^{n-1}}$ and $H\not \cong \Z_2$ as $Q_{2^n}$ has a unique cyclic subgroup of order $2^{n-1}$. The degree of the vertices $\{e\}$ and $\Z_{2^{n-1}}$ is $1$. If $G\cong Q_3$, then it is easy to see degree of $\Z_2$ is $3$, and other vertices are pandent vertices.\\
          
          $\mathbf{G\cong Dic_n}$, where $n=p_1^{a_1}p_2^{a_2}\dots p_k^{a_k}$ and $n\not =2^m$. Then, by definition, dicyclic group of order $4n$ has a unique cyclic subgroup of order $2n$. Let $K$ be the unique cyclic subgroup of order $2n$. Then, every element outside $K$ has order $4$. It is easy to check that $G$ has $n$ cyclic subgroups of order $4$ outside $K$. Also, $G$ contains a cyclic subgroup of order $2p$ for every prime divisor $p$ of $n$. Therefore, the degree of $\Z_2$ is $k+n+1$ as $\{e\}$ is also a maximal subgroup of $\Z_2$. If $n$ is even and $H$ is a subgroup of $K$ of order $p_1^{b_1}p_2^{b_2}\dots p_k^{b_k}\not =2$, then by using the first part the degree of $H$ in $\Gamma (G)$ is $2k-m-r.$ Again, if $n$ is odd and $H$ is a subgroup of $K$ of order $2^{b_{k+1}}p_1^{b_1}p_2^{b_2}\dots  p_k^{b_k}\not =2$, then by using the first part degree of $H$ is $2(k+1)-m-r.$  Moreover, every subgroup of order $4$ outside $K$ is a pendant vertex of $\Gamma (G).$
   \end{proof}

   \begin{cor}\label{cor:maxmin}
       Let $G\cong \Z_n, D_{2n}, Q_{2^n}$ or $Dic_n$, where $n=p_1^{a_1}p_2^{a_2}\dots p_k^{a_k}$, $a_i\geq 1$. Then Table~\ref{tab:MaxMinDegree}, lists the maximum and minimum degree of $\Gamma (G)$.
   \end{cor}
    \begin{table}[ht!]
\centering
\begin{tabular}{|c|c|c|}
 \hline
 $\mathbf{G}$ & $\mathbf{\delta (\Ga (G))}$ & $\mathbf{\Delta (\Ga (G))}$ \\
 \hline
  $\Z_n$ &
   $k$ &  $2k-\ell$, where $\ell$ is the number of\\
              & & \mbox{prime divisors of $n$ with exponent $1$}\\

   \hline
  $D_{2n}$ &
    $1$ & k+n\\
    \hline
  $Q_{2^n}$ & $1$ & $2^{n-2}+1$\\
  \hline
  $Dic_{n}$ & $1$ & $k+n+1$\\

  \hline
 \end{tabular}\\
 \caption{Minimum maximum degree of $\Gamma (G)$}
           \label{tab:MaxMinDegree}
       \end{table}
\begin{proof}
           If $G\cong \Z_n$ and $H$ is a subgroup of $G$ of order $p_1^{b_1}p_2^{b_2}\dots p_k^{b_k}$, then by Theorem~\ref{degree}, the degree of $H$ is minimum if and only if either $b_i=a_i$ or $b_i=0$ for all $1\leq i\leq k$. Thus the minimum degree of $\Gamma (\Z_n)$ is $k$ for all $n$. Again by Theorem~\ref{degree}, $deg (\Z_n)\leq 2k$. Also, if $n$ satisfies the condition $a_i\geq 2$ for all $1\leq i\leq k$, then the subgroup $H$ of order $p_1p_2\dots p_k$ in $\Z_n$ is adjacent to $2k$ subgroups in $\Gamma (\Z_n)$ by using the fact that every maximal subgroup of a cyclic group is of prime index. Therefore the maximum degree of $\Gamma (\Z_n)=2k$, when $a_i\geq 2$ for all $1\leq i\leq k$. Further, let us assume that $\ell$ is the number of prime divisors of $n$ with exponent one, say $p_1,p_2,\dots, p_{\ell}$. If $H$ is any subgroup of $\Z_n$, then there are two possibilities: either $p_i$'s for $1\leq i\leq \ell$ divides $|H|$ or does not divide $|H|$. In both the cases, by Theorem~\ref{degree} degree of $H$ in $\Gamma (\Z_n)$ is $2k-\ell-t$, where $t$ is the number of $p_i$, whose exponent is either zero or $a_i$ in $|H|$ other than $p_1,p_2,\dots, p_{\ell}$. Moreover, if $|H|=p_1p_2\dots p_k$, square-free then by previous theorem $deg (H)=2k-{\ell}$. Hence, in this case, the maximum degree is $2k-{\ell}$, where $\ell$ is the number of prime divisors of $n$ with exponent one.\\
           If $G\cong D_{2n}, Q_{2^n}$ or $Dic_n$, then $\Gamma (G)$ always contains a pendant vertex. Thus, for all these groups, the minimum degree of $\Gamma (G)$ is one. If $G\cong D_{2n}$, then by Table~\ref{tab:degree} the maximum degree of $\Gamma (G)$ is $k+n$ as $n\geq k$. If $G\cong Q_{2^n}$, then by Table~\ref{tab:degree} the maximum degree of $\Gamma (G)$ is $2^{n-2}+1$ because $n\geq 3$. Similarly the maximal degree of $\Gamma (G)$ is $k+n+1$, when $G\cong Dic_n$.
       \end{proof}
       The following corollary is immediate from the proof of Theorem~\ref{degree}.
       \begin{cor}
           Let $G$ be a finite cyclic group and $\Ga (G)$ be the cyclic subgroup graph of $G$ with minimum and maximum degrees $\delta$ and $\Delta$, respectively. Then the degree sequence of $\Ga (G)$ contains all the numbers from $\delta$ to $\Delta$.
       \end{cor}
       \begin{proof}
        Let $G$ be a cyclic group of order $n=p_1^{a_1} p_2^{a_2}\dots p_k^{a_k}$ and $\Gamma (G)$ be the corresponding cyclic subgroup graph. Then using Corollary~\ref{cor:maxmin}, the minimum degree of $\Gamma (G)$ is $\delta (\Gamma (G))=k$ and the maximum degree is $\Delta (\Gamma (G))=2k-\ell$, where $\ell$ is the number of prime divisors of $n$ with exponent $1$. Thus, there are $k-\ell$ prime divisors of $n$ with exponents greater than $1$. Let $H$ be a subgroup of $G$ of order $p_1^{b_1} p_2^{b_2}\dots p_k^{b_k}$ such that $\ell$ prime divisors of $|H|$ have exponents $1$. Then the value of $m(H)=|\{i|b_i=0\}|$ lies in the range $0\leq m(H)\leq k-\ell$. By using Theorem~\ref{degree} $deg_{\Gamma (G)}(H)=2k-m(H)-\ell$. Now it is easy to see that for each value of $m(H)$ from $0$ to $k-\ell$, $deg (H)$ in $\Gamma (G)$ attains all the values from $k$ to $2k-\ell$. Hence, the result holds.
        \end{proof}
        If $n$ is a square-free number, then by Theorem~\ref{degree}, $\Ga(\Z_n)$ is a $k$ regular graph, in fact this graph is well known as a hypercube graph related to boolean algebra. Now, it is interesting to characterize all groups $G$ for which $\Ga(G)$ is a regular graph.
\begin{pro}\label{thm:regular}
    Let $G$ be a finite group. Then $\Gamma (G)$ is regular if and only if $G$ is a cyclic group of square-free order.
\end{pro}
\begin{proof}
If $n$ is a square-free number with $k$ prime divisors, then $\Ga (\Z_n)$ is $k$-regular by Theorem~\ref{degree}.\\
     Let $G$ be a finite group of order $n=p_1^{a_1}p_2^{a_2}\ldots p_k^{a_k}$ such that $\Gamma (G)$ is regular. Then, by Cauchy's theorem, $deg (\{e\})\geq k$. Further, assume that $H$ is a non-trivial cyclic subgroup of $G$ such that $H$ is not contained in any other cyclic subgroup of $G$ and $|H|=p_1^{b_1}p_2^{b_2}\ldots p_k^{b_k}$. Then, by using the fact that every maximal subgroup of $H$ is of prime index and $H$ has a unique subgroup corresponding to every divisor of $|H|$, we have $deg (H)\leq k$. This implies that $deg (\{e\})=k$ and $G$ has a unique subgroup of order $p_i$ for all $1\leq i\leq k$. Thus, every Sylow $p_i$ subgroup of $G$ is either cyclic or generalized quaternion.
 If corresponding to some prime $p_i$, $G$ has more than one Sylow $p_i$-subgroup or order of Sylow $p_i$ subgroup is more than $p_i$, then it is easy to see that either $deg (\{e\})>k$ or $deg (\Z_{p_i})>k$, which is a contradiction. This completes the proof.
\end{proof}
If $G\cong D_{2n}, Q_{2^n}$ or $Dic_n$, where $n\geq 3$, then $\Gamma (G)$ is not Eulerian as it always contains a pendant vertex. By Theorem~\ref{thm:regular}\; if $n$ is square-free and has an even number of prime divisors, then $\Ga(\Z_n)$ is Eulerian. In fact, the converse is also true, that is, if  $\Ga(\Z_n)$ is Eulerian, then $n$ is square-free and has an even number of prime divisors. We will prove this by contradiction. Suppose the exponent is greater than $1$ for at least one prime divisor $p_i$ of $n$, without loss of generality, assume $a_1> 1$. Now consider the subgroups $L$ and $M$ of $\Z_n$, of order $p_1^{a_1}$ and $p_1p_2^{a_2}\dots p_k^{a_k}$ respectively then by Theorem~\ref{degree}, the degree of $L$ is $k$ and degree of $M$ is $k+1$. Therefore, either the degree of $L$ or $M$ is odd, which implies that $\Gamma (\Z_{n})$ is not Eulerian. Hence $\Gamma (\Z_{n})$ is Eulerian if and only if $n$ is square-free and $n$ is the product of an even number of primes.

       The problem of finding all the cyclic subgroups of a finite group is not completely solved. Therefore, it is not easy to calculate the diameter of a cyclic subgroup graph in general. The following result gives bounds on the diameter of $\Ga (G)$ for any finite group $G$.

\begin{theorem}
Let $G$ be a finite group, and let $g$ be an element of $G$ which maximizes the
number of prime divisors of its order (counted with multiplicity); let this
number be $r$. Then $r\le\Ga(G)\le 2r$.
\label{t:diambound}
\end{theorem}

\begin{proof}
Since one step in the cyclic subgroup graph can only change the number of prime
divisors by $1$, it takes at least $r$ steps to reach $\langle g\rangle$ from
the trivial group; but as we have seen, it can be done in $r$ steps. So the
distance from any cyclic subgroup of $G$ to the trivial subgroup is at most $r$,
and this value is achieved. Moreover, given any two cyclic subgroups $H_1$ and
$H_2$, we can move from either one to the trivial subgroup in at most $r$ steps,so we can move from $H_1$ to $H_2$ in at most $2r$ steps. 
\end{proof}

The lower bound is attained by cyclic groups. For that we have seen that $\Ga(\Z_n)$
is isomorphic to a Cartesian product of paths; the distance between two
vertices is the sum of their distances between coordinates in these
paths, which cannot exceed the sum of the lengths of the paths.

The upper bound is attained by $\Z_n\times\Z_n$. To see this, let
$G=A\times B$ where $A\cong B\cong\Z_n$. Let $f(H)$ denote the number of
prime divisors of a cyclic group $H$; and, for a cyclic subgroup $X$ of $G$,
let $m(X)=f(X\cap A)-f(X\cap B)$. Suppose that $X$ and $Y$ are consecutive
subgroups on a shortest path from $A$ to $B$. Then, without loss of
generality, we may suppose that $X\le Y$, with $|Y:X|=p$, where $p$ is a prime
number. Then $|(A\cap Y):(A\cap X)|=1$ or $p$, and similarly for
$|(B\cap Y):B\cap X)|$. Hence $m(X)$ and $m(Y)$ differ by at most~$1$. However,
$m(A)=f(A)$ and $m(B)=-f(A)$; so there must be at least $2r$ steps from $A$
to $B$. Thus, the distance between $A$ and $B$ is at least $2r$, and so the
diameter of $\Ga(G)$ is at least $2r$. Then equality follows from the theorem.

We see that every natural number $r$ is the diameter of the cyclic subgroup
graph of a group: simply take the group $\Z_n$, where $n$ has $r$ prime
factors (counted with multiplicity).

\begin{lemma}
    A cyclic subgroup graph is a cograph if and only if it is a star graph, a path with three vertices, or a cycle of four vertices.
\end{lemma}

\begin{proof}
    If a graph is a star, or a path with three vertices, or a cycle of four vertices, then the corresponding cyclic subgroup graph is a cograph because it is $P_4$ free.\\
    If the graph is not of one of these types, then either it is a tree that is not a star, then it is easy to see that it has $P_4$. Also, if the graph is not a tree, then either it contains a cycle with more than four vertices or a cycle with four vertices and at least one more pendant vertex attached to one of the vertices of the cycle. In this case, also $\Gamma (G)$ contains $P_4$. Hence, the result holds. 
\end{proof}
\section{Cyclic subgroup graph of nilpotent groups}\label{sec:nil}
Now we discuss the cyclic subgroup graph of a nilpotent group. As the class of nilpotent groups has nice properties, it is interesting to examine the properties of the cyclic subgroup graph for this class. One such property of a subgroup graph corresponding to a nilpotent group is that the subgroup graph of a nilpotent group is bipartite: for, if $H<K$, then there is a subgroup $L$ which is between $H$ and $K$ such that $L$ is a maximal subgroup of $K$ and so $K$
has a prime index in $L$. So if we divide the subgroups by the parity of the number of prime divisors, it is a bipartition. Therefore, the subgroup graph of any finite nilpotent group is bipartite. 
\subsection{Cyclic subgroup graph of $p$-group}

Let $G$ be a group of order $p^n$. Then
\begin{enumerate}
    \item Let the exponent of $G$, $exp (G)=p^e$. Then the possible order of a cyclic subgroup of $G$ is $1,p,p^2,\dots, p^e$. The cyclic subgroups in $\Gamma (G)$ can be layered by using their orders as follows: $L_i=\{H\leq G~:~H~\text{is cyclic and}~|H|=p^i\}$, where $0\leq i\leq e$ and $L_0=\{e\}$. In $\Gamma (G)$, edges occur between consecutive layers. 
    \item \textbf{Degree of a vertex}: Let $H$ be a cyclic subgroup of $G$ with $|H|=p^i$. Then\\ $deg_{\Gamma (G)}(H)=$
    $\begin{cases}
            1+|\{K\in C(G):H<K,|K|=p^{i+1}\}| & \mbox{if $i\not =0$},\\
            |\{K\in C(G):H<K,|K|=p^{i+1}\}|  & \mbox{if $i=0$.}\\
              \end{cases}$ 

    Thus, the degree of any maximal cyclic subgroup in $\Gamma (G)$ is $1$. This shows that the minimum degree $\delta (\Gamma (G))=1$.
    \item \textbf{Diameter bound}: The general diameter bound of
Theorem~\ref{t:diambound} holds. But in the case of a nilpotent group, the
element $g$ is the product of elements of maximum order in the Sylow subgroups.
\end{enumerate}

\begin{theorem}\label{thmtree}
Let $G$ be a finite group. Then we have the following:
\begin{enumerate}[label={(\arabic*)}]
\item\label{thm:nil:part1} If $G$ is not a $p$-group, then $\Gamma (G)$ has a pendant vertex if and only if $G$ is a non-nilpotent group.
\item\label{thm:nil:part2} If $G$ is a nilpotent group, then $\Gamma (G)$ is a tree if and only if $G$ is a $p$-group.
\end{enumerate}
\end{theorem}

\begin{proof}
\textbf{\ref{thm:nil:part1}} If $G$ is not a $p$-group and $H$ is a cyclic subgroup of $G$ such that the degree of $H$ in $\Gamma (G)$ is $1$, then by definition of $\Gamma (G)$, $H\cong \Z_{p^k}=\langle a \rangle$, where $k\geq 1$. Since $G$ is not a $p$-group, then there exists a prime number $q\not =p$ such that $G$ has an element of order $q$. Let $b$ be the element of order $q$ in $G$. Moreover, the element $a$ does not commute with the element $b$ as $deg (H)=1$. Hence, $G$ is not nilpotent.\\
Conversely, let $G$ be a nilpotent group, which is not a $p$-group, and let $p$ and $q$ be any two distinct prime divisors of $|G|$. Since $|G|$ has at least $2$ prime divisors, then by Cauchy's theorem $deg (\{e\})>1$. If $H$ is a cyclic subgroup of $G$ such that $|H|$ has at least $2$ prime divisors, then one can easily check that $H$ has at least two maximal subgroups, thus $deg (H)>1$. Also, if $H\cong \Z_{p^k}$, then $H$ is adjacent to subgroups of order $p^{k-1}$ and $p^kq$ as $G$ is nilpotent. Therefore, $G$ has no pendant vertex.\\

\textbf{\ref{thm:nil:part2}} Suppose $G$ is a $p$-group such that $\Gamma (G)$ is not a tree. Then $\Gamma (G)$ contains a cycle say $H_1\sim H_2\sim \cdots \sim H_n \sim H_1$. If $H_i\subseteq H_{i+1}$ for all $1\leq i \leq n$ and $H_n\subseteq H_1$, then $H_n\subseteq H_1 \subseteq H_2\subseteq H_n$, which is a contradiction. Therefore there must exist some $i$,  $1\leq i \leq n$ such that
$H_i\subseteq H_{i+1} \supseteq H_{i+2}$, which is also a contradiction as $H_{i+1}$ is a cyclic group of order $p^k$, where $k\geq 1$ so it contains a unique maximal subgroup. Hence, $\Gamma (G)$ is a tree.\\
Conversely, let $G$ be a nilpotent group, which is not a $p$-group. Then, by using the previous part, $\Gamma (G)$ does not contain any pendant vertex. Hence, the result follows.
\end{proof}
\begin{cor}\label{cor:girth }
   If $G$ is a finite group, then the girth of $\Gamma (G)$ is either $4$ or infinite. 
\end{cor}
\begin{proof}
The proof follows from the proof of the above result and by using the fact that $\Gamma (G)$ is triangle-free.
\end{proof}
   \section{Minimal non-cyclic groups}\label{sec:3}
   Now, we discuss the properties of cyclic subgroup graphs for minimal non-cyclic groups. A non-cyclic group is said to be {\em minimal non-cyclic} if all its proper subgroups are cyclic. These groups are completely characterized in \cite[Proposition~2.8]{jafari2017number}.
\begin{pro}\cite[Proposition~2.8]{jafari2017number}
    Let $G$ be a minimal non-cyclic group. Then $G$ is isomorphic to $\Z_p\times \Z_p$, $Q_8$ or $\langle a,b|a^q=b^{p^r}=1, b^{-1}ab=a^s\rangle$, where $r,s\in \N,q\nmid{s-1}, q|{s^p-1}$ and $p,q$ are distinct primes.
   \end{pro}
   \noindent If $G\cong {\Z_p\times \Z_p}$ or $Q_8$, then it is easy to understand the structure of $\Ga (G)$. Therefore from now onwards we assume that $G=\langle a,b|a^q=b^{p^r}=1, b^{-1}ab=a^s\rangle$, where $r,s\in \N,q\nmid{s-1}, q|{s^p-1}$ and $p,q$ are distinct primes. Also, if $p| (q-1)$, then there always exists a unique minimal non-cyclic group of order $p^rq$. The group $G$ is supersolvable as $G$ contains a cyclic normal subgroup of order $q$ such that $G/H$ is supersolvable.\\
   
   \textbf{Subgroup structure of $G$:}
   The group $G$ has a unique subgroup of order $p^{r-1}q$ by a result given in \cite[Page No.~33]{taunt_1949}. Thus, $G$ has a unique subgroup of order $q$, because if the subgroup of order $q$ is not unique in $G$, then $G$ has a subgroup of order $p^{r-1}q^2$, which is a contradiction. The number of Sylow $p$-subgroups of $G$ is $q$ as the Sylow $q$-subgroup is normal. Also, $G$ has a unique subgroup of order $p^i$ for all $1\leq i \leq{r-1}$  as $G$ has a unique subgroup of order $p^{r-1}q$. Therefore, it is easy to see that all the elements outside the unique subgroup of order $p^{r-1}q$ are of order $p^r$. Therefore the group $G$ contains a unique subgroup corresponding to every positive divisor of $|G|$ except $p^r$ and
 the number of subgroups of $G$ is $2r+q+1$.\\
 
   \noindent The Figure~\ref{fig3} shows the cyclic subgroup graph corresponding to the group $G$.

   \begin{figure}[!htbt]
  \centering
  \begin{tikzpicture}[scale=0.80]
	\vertex (1) at (0,0)[label=below:$\{e\}$]{};
	\vertex (2) at (-2,1)[label=left:$\Z_q$] {};
	\vertex (3) at (2,1) [label=right:$\Z_{p}$] {};
	\vertex (4) at (-2,2)[label=left:$\Z_{pq}$]  {};
	\vertex (5) at (2,2)[label=right:$\Z_{p^2}$]  {};
        \vertex (6) at (-2,3)[label=left:$\Z_{p^2q}$]  {};
	\vertex (7) at (2,3)[label=right:$\Z_{p^3}$]  {};
        \vertex (8) at (-2,4)[label=left:$\Z_{p^3q}$]  {};
         \draw(-2,4.5)node[above]{\Large{$\vdots$}};
        \vertex (9) at (2,5)[label=right:$\Z_{p^{r-1}}$]  {};
        \vertex (10) at (-2,6)[label=left:$\Z_{p^{r-1}q}$]  {};
        \draw(2,3.5)node[above]{\large{$\vdots$}};
        \vertex (11) at (2,6)[label=right:$\Z_{p^r}$]  {};
         \vertex (12) at (1,6)[label=left:$\Z_{p^r}$]  {};
    \draw(3.5,6)node[right]{\Large$\dots$};
          \vertex (13) at (5,6)[label=right:$\Z_{p^r}$]  {};
        
	\path
	  	(1) edge(2)
		(1) edge(3)
		(2) edge(4)
            (3) edge(4)
            (3) edge(5)
		(4) edge(6)
            (5) edge(6)
		(5) edge(7)
		(6) edge(8)
		(7) edge(8)
            (9) edge(10)
            (9) edge(11)
            (9) edge(12)
            (9) edge(13)

	 ;
\end{tikzpicture}
\caption{$\Ga (G)$ for the minimal non-cyclic group of order $p^rq$.}
\label{fig3}
\end{figure}

\vspace{0.2cm}
\noindent \textbf{Properties of $\mathbf{\Ga (G)}$:}
    Let $G$ be a minimal non-cyclic group of order $p^rq$. Then, the corresponding cyclic subgroup graph $\Ga (G)$ has the following properties.
    \begin{enumerate}[label={(\arabic*)}]
    \item The number of vertices in $\Ga (G)$ is $2r+q$.\\

    \item The graph $\Ga (G)$ contains cycles starting at $\{e\}$ of the length $4,6,8,\ldots ,2r$.\\

    \item The diameter of $\Ga (G)$ is $r+1$. Here, we are giving proof for the sake of completeness. If $G$ is a minimal non-cyclic group of order $p^rq$, then the order of any cyclic subgroup of $G$ is either $p^i$ or $p^jq$, where $0\leq i\leq r$ and $0\leq j\leq {r-1}$. Further, $G$ is a CLT group as it is supersolvable, and it contains a unique subgroup corresponding to every divisor of $|G|$ except $p^r$. Since all the subgroups of $p^r$ are maximal, they are all pendant vertices in $\Ga (G)$. If $H$ and $K$ are any two vertices of order $p^i$ and $p^j$, where $0\leq i,j \leq r$, then the distance between them in $\Ga (G)$ is $|i-j|$. Again, if $H$ and $K$ are any two vertices of order $p^iq$ and $p^jq$, where $0\leq i,j \leq {r-1}$, then the distance between them in $\Ga (G)$ is $|i-j|$. Now we are left with the case when $H$ is a vertex of order $p^i$ and $K$ is a vertex of $p^jq$, where $0\leq i \leq r$ and $0\leq j \leq {r-1}$, then the distance between them is $|i-j|+1$. This implies that the maximum distance between any two vertices of $\Ga (G)$ is $r+1$ and is attained when $|H|=p^r$ and $|K|=q$. Hence, the diameter of $\Ga (G)$ is $r+1$.\\

        \item  If $H$ is a proper subgroup of $G$, then the degree of $H$ in $\Ga (G)$ is as follows:\\
        $deg(H)=$ $\begin{cases}
            q+2 & \mbox{if $H\cong \Z_{p^{r-1}}$},\\
             1 & \mbox{if $H\cong \Z_{p^r}$ },\\
              2 & \mbox{if $H\cong  \{e\}, \Z_{q}$ or $\Z_{p^{r-1}q}$},\\
               3 & \mbox{otherwise}.\\
              \end{cases}$\\
    The proof is as follows. A subgroup of order $p^r$ is maximal and cyclic in $G$, and it is adjacent to the subgroup of order $p^{r-1}$ in $\Ga (G)$, which is characteristic in $G$. Therefore, all the subgroups of order $p^r$ are pendant vertices in $\Ga (G)$ as $G$ is not cyclic. Further, the characteristic subgroup of order $p^{r-1}$ is adjacent to $q$ subgroups of order $p^r$ and to a unique subgroup of orders $p^{r-2}$ and $p^{r-1}q$, thus its degree in $\Ga (G)$ is $q+2$. Since $G$ has a unique subgroup of order $p$ and $q$, then the degree of the vertex $\{e\}$ in $\Ga (G)$ is $2$. The subgroup of order $p^{r-1}q$ is maximal in $G$ and it contains two maximal subgroups of index $p$ and $q$, thus its degree in $\Ga (G)$ is $2$. The vertex $\Z_q$ is adjacent to the vertices $\{e\}$ and $\Z_{pq}$, therefore its degree in $\Ga (G)$ is also $2$. If $H$ is a vertex of order $p^iq$, where $1\leq i \leq {r-2}$, then $H$ is adjacent to two of its maximal subgroups of index $p$ and $q$. Also, $H$ is adjacent to the subgroup of order $p^{i+1}q$, therefore degree of $H$ in $\Ga (G)$ is $3$. Similarly if $H$ is a subgroup of order $p^i$, where $1\leq i \leq {r-2}$, then the degree of $H$ in $\Ga (G)$ is also $3$. This completes the proof.\\
    
    \item The degree sequence of $\Ga (G)$ is $\{
   \underbrace{1,1,\ldots ,1,}_\text{q times}2,2,2,
    \underbrace{3,3,\ldots ,3,}_\text{2r-4 times}q+2
   \}$. Hence, $\Ga (G)$ contains $3r+q-2$ edges, and it is neither regular nor Eulerian.\\
   
    \item The independence number of $\Ga (G)$ is $r+q$. Let us prove the argument. By Theorem~\ref{path}, $\Ga (G)$ is bipartite with partite sets $V_1$ and $V_2$ of cardinalty $r+q$ and $r$ respectively. Thus, $\Ga (G)$ has an independent set of cardinality $r+q$. Let $S$ be the largest independent set. Then it contains all the cyclic subgroups of order $p^r$ as all these subgroups are pendant vertices and adjacent to the unique subgroup of order $\Z_{p^{r-1}}$. If we take any two vertices $H$ and $K$ of $S$ not isomorphic to $\Z_{p^r}$, then $|H|\not =|K|$ by using the fact $G$ has a unique subgroup of every order other than $p^r$. This implies that there exists a vertex $L$ in $\Ga (G)$ such that $L$ is adjacent to $H$ but not adjacent to $K$ by using the fact every maximal subgroup of a cyclic group is of prime index, and this is true for any arbitrary subgroups $H$ and $K$ of $S$ of order not equal to $p^,r$. Thus, the cardinality of $S$ is less than or equal to $r+q$. Hence, the independence number of $\Ga (G)$ is $r+q$.\\

    \item The domination number of $\Ga (G)$ is given by\\
    $\gamma (\Ga (G)) =\begin{cases}
            2+\lfloor \frac{2r-3}{4} \rfloor & \mbox{if $r\geq 3$ },\\
             r & \mbox{if $r\leq 2$}.\\
              \end{cases}$\\
    If $r\leq 2$, then it is easy to see that the domination number is $r$. If $r\geq 3$, then by Figure~\ref{fig3} the dominating sets are as follows:\\
     $S=$ $\begin{cases}
            \{\Z_{p^{4l}q},\Z_{p^{4l+2}} |l\in \{0,\ldots,k-1\}\}\cup \{\Z_{p^{r-1}}\} & \mbox{if $r\in \{4k|k\in \N\}$},\\
             \{\Z_{p^{4l}}, \Z_{p^{4l+2}q} |l\in \{0,\ldots,k\}\}\setminus\Z_{4k+2} & \mbox{if $r\in \{4k+1|k\in \N\}$},\\
              \{\Z_{p^{4l}q},\Z_{p^{4l+2}} |l\in \{0,\ldots,k-1\}\}\cup \{\Z_{p^{r-2}q}, \Z_{p^{r-1}}\}  & \mbox{if $r\in \{4k+2|k\in \N\}$},\\
               \{\Z_{p^{4l}q},\Z_{p^{4l+2}} |l\in \{0,\ldots,k\}\} & \mbox{if $r\in \{4k+3|k\in \N\}$}.\\
              \end{cases}$\\
    \end{enumerate}
    \section{From power graph and enhanced power graph to cyclic subgroup graph}\label{sec:relationship}
    Let $G$ be a finite group. The power digraph of $G$ has vertex set $G$ and arcs $x\rightarrow y$ if $y$ is a power of $x$. The power graph of $G$ is the undirected version:
$x\sim y$ is an edge if one of $x$ and $y$ is a power of the other. The enhanced
power graph has an edge $x\sim y$ if there exists an element $z$ such that both $x$ and $y$ are powers of $z$.
\begin{theorem}
    Either the power graph or the enhanced power graph of $G$ determines the cyclic subgroup graph up to isomorphism.
\end{theorem}
\begin{proof}
   We can see that either the power graph or the enhanced power graph determines the directed power graph up to isomorphism by~\cite{cameron2010power,zahirovic2020study,bubboloni2025}. Now it is sufficient to show that the directed power graph determines the cyclic subgroup graph up to isomorphism.

   The relation $x\equiv y$ if $x\rightarrow y$ and $y\rightarrow x$ in the power digraph is an equivalence relation. Moreover, two elements are equivalent if and only if they generate the same cyclic subgroup. So, collapsing each equivalence class to a
single vertex gives the partial order on cyclic subgroups. We note that the edges of the cyclic subgroup
graph are the covering pairs in the partial order. Hence, the result holds.
\end{proof}
\section{Relationship between cyclic subgroup graph and poset}\label{sec:poset}
A partial order on a set is a binary relation $\leq$ which is reflexive,
antisymmetric, and transitive. A partially ordered set or poset is a pair $(X,\leq)$, where $X$ is a set and $\leq$ is a partial order on $X$. As usual, $x<y$ means $x\leq y$ and
$x\not =y$.

A finite poset $(X,\leq)$ is ranked if, for every $x\in X$, all maximal descending chains $x=x_0>x_1>\dots>x_r$ have the same length $r=r(x)$, which is the
rank of $x$.

A poset may have two important additional properties:
\begin{enumerate}
\item\label{item(1)} It is ranked and has a unique minimal element.
\item\label{item(2)} It is a meet-semilattice, this means that any two elements $x$ and $y$ have a unique meet or greatest lower bound. There is an element $z\in X$ is such that $z\leq x$ and $z\leq y$, and if $w\leq x$ and $w\leq y$ then $w\leq z$.
\end{enumerate}
The comparability graph of a poset $(X,\leq)$ is the graph with vertex set $X$ and edges $x$ is connected to $y$ whenever $x<y$ or $y<x$. Comparability graphs of posets have several special properties. For example, they are perfect, which means that every induced subgraph has the clique number equal to the chromatic number. 

In a poset $(X,\leq)$, we say that $x$ covers $y$ if $y<x$ but there is no element $z$ with $y<z<x$. The covering graph of the poset $(X,\leq)$ has vertex set $X$, with $x$ is adjacent to $y$ if $x$ covers $y$ or $y$ covers $x$. It is a spanning subgraph of the comparability graph.

If we know the orders of covering pairs in a finite poset, then the entire poset is determined: for $x\leq y$ if and only if there is a chain
$x=x_0\leq x_1\leq\dots \leq x_r=y$
in which $x_i$ covers $x_{i-1}$ for $1\leq i\leq r$.

The problem of deciding which graphs are covering graphs of posets is
known to be NP-hard~\cite{brightwell1993complexity}. However, for posets satisfying condition~\ref{item(1)} above, the problem is much easier.
\begin{pro}\label{poset(1)}
    If $(X,\leq)$ is a finite ranked poset with a unique minimal
element, then the covering graph is bipartite and connected.
\end{pro}
\begin{proof}
    The covering graph is connected since every element has a path to the unique minimal element. It is bipartite since, in any covering
pair, the ranks of the two elements differ by $1$, so if we take two sets of vertices based on the parity of the rank, it defines the bipartition.
\end{proof}
\begin{pro}\label{poset(2)}
    If $B$ is a connected bipartite graph with vertex set $X$, and $v$ is any vertex of $B$, then there is a ranked poset $(X,\leq)$ with least element $v$ and covering graph $B$.
\end{pro}
\begin{proof}
    We define the rank function and the elements covered by a given
element by induction. The induction starts with the only element of rank zero, which is $v$, and nothing is covered by $v$. Suppose that the elements of rank $r$ and the elements they cover have been determined. Choose an element $x$ of rank
$r$. The elements covered by $x$ are known; the remaining neighbours of $x$ have rank $r+1$ and cover $x$.
\end{proof}
A Hasse diagram of a finite poset is a drawing of the covering graph in
the plane such that, if $y$ covers $x$, then the position of $y$ is higher up than that of $x$ (that is, its $Y$-coordinate is greater than that of $X$). Clearly, the Hasse diagram determines the order of all covering pairs and hence the order
of the whole poset.

\textbf{Cyclic subgroup graph as a poset~:}
Let $G$ be a finite group. The cyclic subgroup poset of $G$ has as elements of the set $C(G)$ the collection of cyclic subgroups of $G$, and the partial order is the subgroup relation. The cyclic subgroup graph is the covering graph of this poset.
\begin{ex}\label{example2}
    Consider the group $Q_8=\{1, -1, i, -i, j, -j, k, -k\}.$ The set of all cyclic subgroups of $G$ is $C(G)=\{\{e\}, \Z_2=\langle -1\rangle, \Z_4=\langle i\rangle, \Z_4=\langle j\rangle, \Z_4=\langle k\rangle\}$, and the corresponding cyclic subgroup graph $\Gamma (G)$ is given in Figure~\ref{fig12}.
    \begin{figure}
    \begin{center}
\begin{tikzpicture}
[
    every node/.style={
        draw,
        rounded corners,
        inner sep=5pt,
        font=\large
    },
    every edge/.style={
        draw,
        ->,
        >=stealth,
        thick
    }
]

\node (E)  at (0,0)   {$\{e\}$};
\node (Z)  at (0,1.5) {$\langle -1\rangle$};

\node (I)  at (-3,3)  {$\langle i\rangle$};
\node (J)  at (0,3)   {$\langle j\rangle$};
\node (K)  at (3,3)   {$\langle k\rangle$};

\path
(E) edge (Z)

(Z) edge (I)
(Z) edge (J)
(Z) edge (K);

\end{tikzpicture}

\caption{Corresponding Graph of $Q_8$}
\label{fig:Q8-hasse-diagram}
\end{center}
\label{fig12}   
\end{figure}
\end{ex}

Now it is easy to see that the cyclic subgroup poset is bipartite and connected.
\begin{pro}
    The cyclic subgroup poset is bipartite and connected.
\end{pro}
\begin{proof}
    The proof is straightforward, as the cyclic subgroup poset satisfies properties~\ref{item(1)} and~\ref{item(2)}, that is, it is a ranked poset with a unique minimal element, and is a meet semilattice. The minimal element is the trivial subgroup $\{e\}$ and
the rank of a cyclic subgroup is the number of prime divisors of its order
(counted with multiplicity). The meet of two cyclic subgroups $H$ and $K$ is their intersection $H\cap K$. Therefore, by using Proposition~\ref{poset(1)}, the cyclic subgroup poset is bipartite and connected.
\end{proof}
In the context of a cyclic subgroup poset, it is interesting to ask the following problem.
\begin{problem}
    When is it the case that different choices of the orbits containing the trivial subgroup give rise to posets of cyclic subgroups of groups?
\end{problem}
We do not know the complete answer to the above problem. The following Lemma gives the answer to the problem for the family of star graph.
\begin{lemma}
    Let $\Gamma (G)=K_{1,n}$ be a star graph if the orbit of the identity vertex has been changed and the graph structure is preserved. Then $n=4$ and $\Gamma (G)=K_{1,4}$.
\end{lemma}
\begin{proof}
   Let $G$ be a group and $\Gamma (G)$ be a star graph $K_{1,n}$. First, suppose that the identity $\{e\}$ is a leaf in $\Gamma (G)$. Then the degree of identity $\{e\}$ in $\Gamma (G)$ is one. Also, in $\Gamma (G)$ the neighbours of $\{e\}$ are subgroups of prime order. This implies that $G$ has a unique subgroup of prime order. Let $H\cong \Z_p$ be the unique subgroup of $G$ of prime order. Also $H$ is the center of the star graph so all other cyclic subgroups of $G$ are connected to $H$. Let the other leaves except identity in $\Gamma (G)$ are $B_1, B_2,\dots B_{n-1}$. Then $H\leq B_i$, where $1\leq i\leq n-1$. Also, $H$ is the proper maximal subgroup of each $B_i$.

   If at least two primes divide the order of $B_i$, then $\Gamma (G)$ contains a cycle. Thus $|B_i|=p^2$ for all $1\leq i\leq (n-1)$. Since $G$ has a unique subgroup of prime order so either $G$ is either cyclic or generalized quaternion. If $G$ is cyclic then $G\cong \Z_{p^2}$ and $\Gamma (G)=K_{1,2}$. But in this case, the vertex $\{e\}$ can not occur as a central vertex, so this case is not possible. Hence $G\cong Q_8$ and $\Gamma (G)=K_{1,4}$.
   
   If $\Gamma (G)$ is $K_{1,4}$ and identity $\{e\}$ is the central vertex, then $G$ has exactly four subgroups of prime order. By Sylow theorem either, all the non-trivial cyclic subgroups of $G$ are of order $3$ or three subgroups are of order $2$ and one subgroup is of order $3$. Hence $G\cong \Z_3\times \Z_3$ or $S_3$ symmetric group of order $3$.  
\end{proof} 
   \section{Acknowledgement}
    The authors wish to thank Dr. Angsuman Das for the helpful discussions and suggestions during this work. The third-named author is supported by the University Grant Commission (UGC), India, under the scheme UGC-SRF.\\

\noindent \textbf{Conflict of interest:} The authors declare that they have no known competing financial interests or personal
relationships that could have appeared to influence the work reported in this paper.\\

\bibliographystyle{plain}
  \bibliography{refs}
    \end{document}